\documentclass[final]{siamltex}

\def\bx{{\bf x}}

\def\C{{C\kern-.647em I}}

\def\bv{{\mathbf v}}
\def\bx{{\mathbf x}}
\def\mod{{\rm mod}}

\def\h{{\ \buildrel \rm h \over = \ }}

\def\R{I\!\!R}

\def\beq{\begin{equation}}
\def\eeq{\end{equation}}

\usepackage[dvips]{graphicx}
\begin{document}
\title{The Spectral Basis and Rational Interpolation\thanks{This 
        work was supported by NIP of UDLA. The author is a member of Sistema
Nacional de Investigadores, Expiente No. 14587.}}
\author{Garret Sobczyk \thanks{Departamento de F\'isica y Matem\'aticas,
Universidad de las \'Americas - Puebla, Mexico,
 72820 Cholula, M\'exico, (garrete.sobczyk@udlap.mx).}}
\maketitle
 \begin{abstract} The Euclidean Algorithm is the often forgotten key to rational approximation techniques,
 including Taylor, Lagrange, Hermite, osculating, cubic spline, Chebyshev, Pad\'e and other interpolation schemes.
  A unified view of  these various interpolation techniques is eloquently expressed in terms of the 
  concept of the spectral basis of a factor ring of polynomials. When these methods are applied to the minimal polynomial
  of a matrix, they give a family of rational forms of functions of that matrix.
  \end{abstract}
\pagestyle{myheadings}
\thispagestyle{plain}
\markboth{G. Sobczyk}{RATIONAL INTERPOLATION}
\section{The Euclidean Algorithm and Spectral Basis}
   The euclidean algorithm has many important well-known consequences in number theory, algebra and analysis.
In the spirit of \cite{AMS}, we are mainly interested here in some of its consequences regarding interpolation of functions over the real
or complex numbers. Let $\R[x]$ and $\C[z]$ denote the rings of real-valued and complex-valued polynomials over
the real and complex number fields $\R$ and $\C$, respectively. Whereas we state our results in terms of polynomials
over the field of real numbers $\R$, all of the results are equally valid for polynomials over $\C$.

Let $h(x)$ denote the monic polynomial defined by
   \beq h(x)=\prod_{i=1}^r(x-x_i)^{m_i}, \label{divisor} \eeq
where $\{x_1, \ldots, x_r\}$ are the distinct real roots of $h(x)$ with multiplicities $\{m_1, \ldots, m_r\} $, 
respectively. Let 
$f(x)\in \R[x]$ be any polynomial in $\R[x]$. Then the euclidean algorithm simply tells us that there will
always exist polynomials $g(x)$ and a remainder $r(x)$ such that
   \beq  f(x) = g(x)h(x)+r(x) \label{EA} \eeq
where $0\le \deg(r(x))< \deg(h(x)) $ or $r(x) \equiv 0$. When equation (\ref{EA}) holds, we say that
$f(x)=r(x)\ \mod(h)$ where $h=h(x)$, or more concisely, that
$f(x)\h r(x)$. We denote the ring of all real polynomials modulo $h(x)$ by $\R[x]_{h}$.
In the terminology of factor rings, $\R[x]_{h}\widetilde=\R[x]/\!\!<\!h(x)\!>$, meaning that $\R[x]_{h}$ is 
isomorphic to the factor ring
$\R[x]/\!\!<\!h(x)\!>$ of $\R[x]$ generated by the principal ideal $<h(x)>$,  \cite[p.266]{Gal}. Thus, $\R[x]_{h}$ has the
structure of a ring with addition and multiplication of polynomials defined modulo $h(x)$.

   By the {\it standard basis} of {\it residue classes} of $\R[x]_{h}$ we mean
   \beq  {\cal B}_h=\{1,x,x^2,\ldots, x^{m-1} \}, \label{b} \eeq
where $m=\deg(h)=m_1+\ldots +m_r $. However, calculations in  $\R[x]_{h}$ are much more simply carried out 
by appealing to the special properties of the {\it spectral basis} \cite{S0}, \cite{S3}. 
The {\it spectral basis} ${\cal S}_h$ of $\R[x]_{h}$
consists of {\it idempotents} $s_i=s_i(x)$, and {\it nilpotents} $q_i=q_i(x)$ and their powers $q_i^k=q_i^k(x)$,
\beq   {\cal S}_h=\{s_1,q_1,\ldots, q_1^{m_1-1},s_2,q_2,\ldots, q_2^{m_2-1}, \ldots ,s_r,q_r,\ldots, q_r^{m_r-1} \}, \label{sb} \eeq
which satisfy the following properties under addition and multiplication in $\R[x]_{h}$ modulo $h(x)$:

\begin{remunerate} 

\item[\bf Property 1.]  $ s_1 + s_2 + \cdots + s_r = 1 $, and $s_i s_j \h \delta_{ij} s_i$ for $i,j=1, \ldots , r$  where
$\delta_{ij}=0$ for $i\ne j$ and $\delta_{ij}=1$ for $i=j$.

\item[\bf Property 2.]  $q_i s_i \h  q_i$, and $q_i^{m_i-1} \ne 0 \ \mod(h)$ but $q_i^{m_i} \h 0$, for $i=1, \ldots, r$.

\item[\bf Property 3.]  For each $f(x) \in \R[x]$, $f(x) s_i \h \big( f(x) \mod[(x-x_i)^{m_i}]\big)s_i$ for $i=1, \ldots, r$.

\end{remunerate}

    Property 1, shows that the $s_i(x)$ are mutually annihilating idempotents which partition unity. Property 2, shows 
that $s_i$ acts as an identity element when multiplied by the nilpotent $q_i(x)$, and that $q_i(x)$ is a nilpotent
of index $m_i$ in  $\R[x]_h$. Property 3, shows that for each polynomial $f(x)\in \R[x]$, $s_i(x)$ acts as the
 projection of
$f(x)$ onto the ring of polynomials $\R[x]_{(x-x_i)^{m_i}}$ modulo $(x-x_i)^{m_i}$.  
  
   From these three algebraic properties,
we can explicitly {\it solve} for the polynomials that make up the spectral basis. 
   For each $i=1, \ldots, r$, define $h_i=h_i(x)=h(x)/(x-x_i)^{m_i}$. Using Properties 1 and 3, we find that
  \beq  h_i s_i \h h_i ,  \label{solve}  \eeq
 and since $h_i(x_i) \ne 0$, it follows that $h_i^{-1} \ne 0 \ \mod \ (x-x_i)^{m_i}$. Multiplying both sides
of equation (\ref{solve}) by $h_i^{-1} \mod\ (x-x_i)^{m_i}$, we find that 
  \beq s_i(x)=(h_i^{-1}\mod \ (x-x_i)^{m_i})h_i(x), \label{idempotents} \eeq
which gives an explicit solution by taking the first $m_i$ terms of the Taylor series for $h_i^{-1}$ around the point
 $x=x_i$. Having found the idempotents $s_i$, the corresponding nilpotents $q_i$, and their powers, are specified by
   \[ q_i^k  :\h (x-x_i)^k s_i = (x-x_i)^k  (h_i^{-1}\mod \ (x-x_i)^{m_i})h_i(x) \]
      \beq          \h       (h_i^{-1}\mod \ (x-x_i)^{m_i-k})h_i(x), \label{nilpotents} \eeq 
  for $k=0,1, \ldots, m_i-1$. Note that for $k=0$, we get the correct convention that $q_i^0 \h s_i$. 

     The transition from the standard basis (\ref{b}) to the spectral basis (\ref{sb}) is accomplished by first noting
that
   \beq x = \sum_{i=1}^r xs_i = \sum_{i=1}^r ((x-x_i)+x_i)s_i \h  \sum_{i=1}^r (x_i+q_i)s_i , \label{x} \eeq
from which it follows that  
 \[ x^k \h \sum_{i=1}^r (x_i+q_i)^k s_i   \h  \sum_{i=1}^r \sum_{j=0}^{m_i-1} \pmatrix{k \cr j} x^{k-j}_i q_i^j  ,  \]
for $k=0,1,\ldots, m-1$, as easily follows by referring to the properties of the spectral basis \cite{S3}.

 \section{Rational Interpolation}

   Let $f(x)$ be a real-valued function which is continuous and has derivatives to the 
orders  $\{m_1-1, \ldots, m_r-1\} $ at the respective points $\{x_1, \ldots, x_r\}$,
where $\{x_1, \ldots, x_r\}$ are the distinct real roots of $h(x)$ with multiplicities $\{m_1, \ldots, m_r\} $
 as was defined in (\ref{divisor}) of the previous section. 

The function $f(x)$ of the real variable $x$ can be extended to a function of the variable $x\in \R[x]_h$ 
by simply substituting (\ref{x}) into $f(x)$, getting
   \beq  f(x):=f(\sum_{i=1}^r (x_i+q_i)s_i)  \h \sum_{i=1}^r f(x_i+q_i)s_i. \label{fmodh1} \eeq
If we now expand $f(x_i+q_i)$ in a Taylor series about $x=x_i$, we get the desired expression
   \beq  f(x_i+q_i) \h \sum_{k=0}^{m_i-1} \frac{1}{k!}f^{(k)}(x_i)q_i^k , \label{fmodh2} \eeq
where as usual, $f^{(k)}(x_i)=\frac{d^k}{dx^k}f(x)|_{x=x_i}$. Although we have {\it derived}
equations (\ref{fmodh1}) and (\ref{fmodh2}) modulo $h(x)$ from the basic properties of the spectral basis (\ref{sb}),
we could equally well have taken (\ref{fmodh1}) and (\ref{fmodh2}) to be the {\it definition} of $f(x)$ modulo $h(x)$. 

    The interpolation polynomial
    \beq g(x)= f(x) \ \mod\ h(x) \h  \sum_{i=1}^r [\sum_{k=0}^{m_i-1} \frac{1}{k!}f^{(k)}(x_i)q_i^k ]s_i 
                                  \label{osculating} \eeq 
is called the {\it Birkhoff} or {\it osculating} interpolation polynomial of $f(x)$ with respect to $h(x)$. 
 In the special case that $m_1=\cdots = m_r =1$, $g(x)$
is called the {\it Lagrange} interpolation polynomial of $f(x)$, and in the special case when
$m_1= \ldots = m_r = 2$, $g(x)$ is called the {\it Hermite} interpolation polynomial of $f(x)$, \cite[pps.278,287]{Ham},
\cite[p.52]{Sto}. 
When $r=1$,
$g(x)$ reduces to the first $m_1-1$ terms of the Taylor series of $f(x)$ about $x=x_1$.

      Rational interpolation also takes an equally eloquent form when expressed in terms of the spectral basis 
\cite[p.58]{Sto},\cite{Davis}.  Let
$a(x)=\sum_{i=0}^{m-1} a_i x^i$, and $b(x)=\sum_{j=0}^{m-1} b_j x^j$ be polynomials over the real numbers $\R$.
We say that  
      \beq g(x)=\frac{a(x)}{b(x)} \label{rat} \eeq
is a {\it rational} interpolate of $f(x)$ at the points (nodes) $\{x_1, \ldots, x_r\}$  
with multiplicities $\{m_1, \ldots, m_r\}$ if 
   \beq f(x)b(x)-a(x) \h 0.  \label{rat1} \eeq
The usual way of defining rational interpolation involves the solution of a system of linear equations
for coefficients of the polynomials $a(x)$ and $b(x)$. Our definition is simpler
and more direct in that it only requires that the modular relation (\ref{rat1}) holds. When
$b(x)$ has no common zeros with $h(x)$, the equation (\ref{rat1}) is equivalent to 
   \beq f(x)\h a(x)b(x)^{-1}. \label{rat2} \eeq

 Essentially, each choice of the {\it shape polynomial} $b(x)$ in (\ref{rat1}) and (\ref{rat2}) determines a different rational interpolation $g(x)$ 
of $f(x)$, \cite{QW}.  
Because of the homogeneous nature of the rational interpolate (\ref{rat}), we can require that
$b(x_e)=1$ at some point $x_e$ which is chosen not to be one of the roots $x_i$ of $h(x)$ nor a
zero of $f(x)$.   
 Whereas equation (\ref{rat1}) is
defined modulo($h(x)$) using (\ref{fmodh1}) and (\ref{fmodh2}), equation (\ref{rat}) is an {\it ordinary} equality.

    Chebyshev and other kinds of rational interpolation are defined simply by replacing the powers of $x^k$ of
 $x$ in $a(x)$ and $b(x)$ by the corresponding Chebyshev or other sets of orthogonal polynomials of the same order.
 A comprehensive study of the algebraic structure of rational functions has been undertaken
by Luis Verde-Star in a series of papers \cite{V97a, V97b, V97c}. 
Examples of the various kinds of interpolation will be given in Section 3.
 
  Cubic spline interpolation can be defined parametrically in terms of the spectral basis ${\cal S}_{2,2}$
of $\R[t]_{h}$ for $h=h(t)=t^2(t-1)^2$. Using the formulas (\ref{idempotents}) and (\ref{nilpotents})
 from the previous section,
we find that  
   \beq {\cal S}_{2,2}=\{s_1=(2t+1)(t-1)^2,q_1=t(t-1)^2,s_2=(3-2t)t^2,q_2=t^2(t-1) \}. \label{spectral22} \eeq
The piecewise {\it natural cubic spline} $\{g_1(t_1),g_2(t_2),\ldots , g_{k-1}(t_{k-1}) \} $, for $0\le t_i < 1$ and $k\ge 3$,
connecting the successive points $\{ \bx_1, \bx_2, \ldots , \bx_{k} \}$ in $\R^n$, is
defined by
  \beq g_i(t_i)=\bx_i s_1(t_i) + \bv_{i} q_1(t_i)+  \bx_{i+1} s_2(t_i) + \bv_{i+1} q_2(t_i),  \label{cubicspline} \eeq
with the requirements that  
   \beq g_1^{\prime \prime}(0)=0= g_{k-1}^{\prime \prime}(1) \ \ {\rm and} \ \  g_i^{\prime \prime}(1)= g_{i+1}^{\prime \prime}(0)    
    \label{free} \eeq
 for $i=1, \ldots ,k-2$. Taking the second derivatives
 of $g_i(t_i)$, and evaluating at $t_i=0,1$ gives
   \beq  g_i^{\prime \prime}(0) = 6(\bx_{i+1}-\bx_i)-4\bv_i - 2 \bv_{i+1},  \label{second0} \eeq  
and   
    \beq      g_i^{\prime \prime}(1) = -6(\bx_{i+1}-\bx_i)+2\bv_i + 4 \bv_{i+1}. \label{second1} \eeq
The resulting $k$ linear vector equations are uniquely solved for the $k$-unknown tangent vectors $\bv_1, \bv_2, \ldots, \bv_k$. 

If, instead of the requirement (\ref{free}), the vectors $\bv_1$ and $\bv_{k}$ are taken as given, 
and the remaining $(k-2)$ linear vector equations
 \beq  g_i^{\prime \prime}(1)= g_{i+1}^{\prime \prime}(0),    
    \label{bounded} \eeq
for $i=1, \ldots ,k-2$, are uniquely solved for the $(k-2)$ unknown tangent vectors $\bv_2, \bv_2, \ldots, \bv_{k-1}$, the
resulting solution is called the {\it bounded cubic spline}, \cite[p.93]{Sto}.          
 
     Various kinds of {\it rational cubic splines} can also be easily constructed by replacing the spectral basis 
  ${\cal S}_{2,2}$ in (\ref{cubicspline}) by a {\it rational spectral basis} of the form
     \beq {\cal R}_{2,2}(b)=\{s_{r1},q_{r1},s_{r2},q_{r2} \} 
                        \label{ratspectral} \eeq
  for
     \[  s_{r1}(t)=\frac{b(t) s_1(t) \ \mod(h(t))}{b(t)} , \  q_{r1}(t)=\frac{b(t) q_1(t) \ \mod(h(t))}{b(t)}  \] 
  and     
     \[  s_{r2}(t)=\frac{b(t) s_2(t) \ \mod(h(t))}{b(t)} , \  q_{r2}(t)=\frac{b(t) q_2(t) \ \mod(h(t))}{b(t)},  \] 
  where $b=b(t)=1+b_1 t+b_2 t^2 + b_3 t^3$ and $b(1)\ne 0$. Clearly, the rational spectral basis 
  ${\cal R}_{2,2}(b)$ reduces to the ordinary spectral basis  ${\cal S}_{2,2}$, given in (\ref{spectral22}),
 for $b=b(t)=1$. Of course, when
  using the rational spectral basis, the second derivatives $g_i^{\prime \prime}(0)$ and $g_i^{\prime \prime}(1)$,
  given in (\ref{second0}) and (\ref{second1}), must be recalculated. 
  
  \section{Circles}
  
      A circle and other conics are good geometric figures on which to carry out interpolation experiments. 
We derive here several approximations for the unit circle, centered at the origin, using rational spectral bases.

    The rational spectral basis ${\cal S}_h=\{s_{r1},s_{r2},s_{r3} \}$ for $b=1+b_1 t+b_2 t^2$ and $h(t)=(t+1)t(t-1)$,
is defined by
    \[ s_{r1}=\frac{\frac{1}{2}(1-b_1+b_2)t(t-1)}{1+b_1 t+b_2 t^2} ,s_{r2}=-\frac{(t-1)(t+1)}{1+b_1 t+b_2 t^2}  , \]
and
    \[         s_{r3}= \frac{\frac{1}{2}(1+b_1+b_2)t(t+1)}{1+b_1 t+b_2 t^2}.    \]
We wish to optimize (in the sense of {\it least squares}) the choice of $b_1$ and $b_2$
 so that the interpolating curve $g(t)=f(-1)s_{r1}+f(0)s_{r2}+f(1)s_{r3}$ 
to the semi-circle 
$f(t)=\pmatrix{\cos(\pi t/2),\sin(\pi t/2)}$, for $-1 \le t \le 1$ is as good as possible. The values of
$b_1$ and $b_2$ can easily be found by requiring that $g(1/2)\cdot g(1/2) =1=g(-1/2)\cdot g(-1/2)$, giving
the values $b_1=0$ and $b_2=\pm 1$. The value $b_2=-1$, gives the single point $(1,0)$, whereas $b_2=1$
gives the well-known parameterization of the circle $g(t)=(\frac{1-t^2}{1+t^2}, \frac{2t}{1+t^2})$.

    A family of rational approximations to the quarter unit circle through the nodes $(1,0)$ and $(0,1)$, 
with the initial and terminal tangent velocity vectors $(0,v)$ and $(-v,0)$, is specified by
   \[ g(t)=(1,0)s_{r1}+(0,v)q_{r1}+(0,1)s_{r2}+(-v,0)q_{r2}, \] 
 in the rational spectral basis (\ref{ratspectral}). One popular construction of the circle is based on 
NURBS (nonuniform rational B-splines) \cite[p.110]{Sto}. 
 Letting $b_1=-2+\sqrt{2}=-b_2$, the
 choice $v=\sqrt{2}$ precisely eliminates the $t^3$ term in the numerator, and gives the nurb parameterization
   \[ g(t)=(1,0)s_{r1}+(0,\sqrt{2})q_{r1}+(0,1)s_{r2}+(-\sqrt{2},0)q_{r2} \]
    \[ =\Big(\frac{1+(-2+\sqrt{2}) t +(1-\sqrt{2}) t^2}{ 1+(-2+\sqrt{2}) t +(2-\sqrt{2} ) t^2}, 
    \frac{\sqrt{2} t +(1-\sqrt{2}) t^2}{ 1+(-2+\sqrt{2}) t +(2-\sqrt{2} t^2}\Big). \] 
    
    A quite different perfect parameterization of the unit circle through the interpolation points
    $(1,0)$ and $(0,1)$, and taking the initial tangent vector at $(1,0)$ to be $(0,\pi/2)$, can be derived by using
the rational spectral basis of
  \[{\cal S}_{2,1}=\{s_1=-(t+1) (t-1),q_1=-(-1 + t) t,s_2=  t^2   \}    \] 
for $h(t)=t^2(t-1)$. We find that
   \[ g(t)=\frac{(1,0)s_1+\Big(b_1(1,0)+(0,\frac{\pi}{2})\Big)q_1+(0,1)(1+b_1+b_2)s_2}{1+b_1 t+b_2 t^2} \]
Optimizing $b_1$ and $b_2$, we find that
    \[ g(t)=\Big( \frac{8+4(\pi -4)t-4(\pi-2)t^2}{8+4(\pi -4)t+(\pi^2-4 \pi +8)t^2},  
                    \frac{4\pi t+(\pi-4)\pi t^2}{8+4(\pi -4)t+(\pi^2-4 \pi +8)t^2} \Big). \]  
  
      Whereas the above parameterizations give perfect circles, there are many other parameterizations that are
  interesting. For example, consider the family of approximations to the unit semicircle through the points 
$(1,0)$ and $(-1,0)$, given by
     \[ g(t)=(1,0) s_{r1}+(0,v)q_{r1} + (-1,0) s_{r2}+(0,-v)q_{r1} \] 
 employing the rational spectral basis of the kind (\ref{ratspectral}). Choosing $b=1-t+t^2$, and $v=3$
   gives a very good approximation to the unit semicircle for 
$0\le t \le 1$,   \beq g(t) = \Big( \frac{1 - t - 3 t^2+2 t^3}{1-t+t^2},\frac{-3(t-1)t}{1 - t +  t^2}\Big)\widetilde =
        (\cos\pi t,\sin\pi t). \label{cossin} \eeq
  with a least square error less than $.000071$, see figure 3.1 
   \begin{figure}
\begin{center}
\includegraphics[scale=.35]{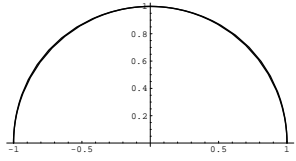}
\caption{The unit semicircle is shown together with its approximation.}
\end{center}
\end{figure} 
   \begin{figure}
\begin{center}
\includegraphics[scale=.5]{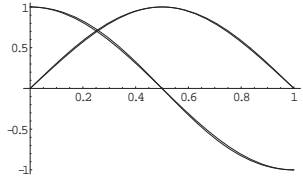}
\caption{Both sine and cosine curves are shown together with their approximations.}
\end{center}
\end{figure}
It is interesting to note
   that this parameterization also gives a good approximation to $\cos \pi t$ and $\sin \pi t$ for
$0\le t \le 1$, see figure 3.2
The series expansions for the approximations to $\cos \pi t$ and
   $\sin \pi t$ are
    \[   \cos \pi t \widetilde = 1-4t^2-2t^3+\sum_{k=1}^\infty (-1)^{k+1}[2t^{3k+1}+4t^{3k+2}+2t^{ 3k+3}] \]
and                        
   \[   \sin \pi t \widetilde = 3t+3\sum_{k=1}^\infty  (-1)^{k}[t^{3k}+t^{3k+1}], \]
which are interesting in their own right.   
   \section{Matrices}
       Let $A$ be any $n\times n$ matrix over a field $\cal K$. The field $\cal K$ may be the real or complex numbers, 
 or even a finite Galois field. It is well known that every matrix satisfies it's characteristic polynomial, defined by
   \[ \varphi(x)=\det (x I - A), \]
   where $I$ is the identity $n\times n$ matrix \cite{HJ90,Gant}. The theory of a spectral basis is directly applicable to a matrix
   whenever the characteristic polynomial is of the form $\varphi(x)=\prod_{i=1}^r(x-x_i)^{n_i}$ for distinct
   roots $x_i \in {\cal K}$. When applying the spectral basis to a matrix $A$, it is better to use the
   {\it minimal polynomial} $\psi(x)=\prod_{i=1}^r(x-x_i)^{m_i}$, where $1\le m_i \le n_i$ for each $i=1, \ldots , r$.
   The minimal polynomial of the matrix $A$ is defined by the condition that it is unique monic polynomial of least degree
   for which $\psi(A)=0$.
   
     Any of the interpolation formulas, developed in terms of the spectral basis in the previous sections,
 apply immediately to the matrix $A$, provided that $h(x)=\psi(x)$. This is because the relationship that
$h(x) = 0 \ \mod \ h$ is precisely reflected in the condition that $\psi(A)=0$ for the minimal
polynomial $\psi$ of the matrix $A$. 
Thus, the {\it spectral form} of the matrix $A$ is given by simply replacing $x$ by
   the matrix $A$ in (\ref{x}), getting  
   \beq A =   \sum_{i=1}^r x_i s_i(A)+q_i(A)=\sum_{i=1}^r x_i S_i+Q_i, \label{A} \eeq
 for $S_i=s_i(A)$ and $Q_i=q_i(A)$. The matrices $\{S_i,Q_i \}$ of the spectral basis satisfy exactly
 the same rules under matrix addition and multiplication as does the polynomials $\{s_i,q_i\}$ of
 the spectral basis modulo $h(x)$, \cite{S0}.
 
 For example, the matrix
     \[  A = \pmatrix{-2 & -1 & 0 & 0 \cr -1 & 5 & 1 & 1 \cr -6 & -5 & -1 & 0 \cr -8 & -10 & -3 & 0} \]
  has both characteristic and minimal polynomials $h=t^2(t-1)$. As a consequence, we can apply all the
  above interpolation formulas, found for the rational spectral basis of $h(x)$, without modification. Using
(\ref{osculating}), (\ref{spectral22}), (\ref{cossin}) and (\ref{A}), we find that
  \[ \cos{\pi A} = I - 6 A^2+4 A^3 = \frac{I - A-3 A^2+2 A^3}{I-A+A^2} =
      \pmatrix{1& 2& -10& 6\cr -10& -3& 26& -14\cr 10 & 8& -45& 26\cr 20& 14& -82& 47}\]
and  
  \[ \sin{\pi A} =  A-A^2=\frac{-3A^2 +3A}{I-A+A^2}=
\pmatrix{21& -6& -3& -3\cr -48& 18& 0& 12\cr 87& -27& -9& -15\cr 156& -51& -12& -30}.\]
As a check, we calculate $\cos^2 \pi A + \sin^2 \pi A = I$ as expected.

\section*{Acknowledgments}
The author gratefully acknowledges the support of Dr. 
Gerardo Ayala of NIP 
and Dr. Andres Ramos of the Department of Mathematics at UDLA. He also thanks
his student, Omar Le\'on S\'anchez, for his help on the approximations
to the circle.

\end{document}